\documentclass[12pt]{article}
\usepackage{amssymb,amsfonts,euscript,mathrsfs}
\usepackage{amsmath}
\usepackage[cp1251]{inputenc}
\usepackage[russian] {babel}

\begin{document}
\large
\begin{center}
 \bf V.N. Dubinin\\
 ON THE LEMNISCATE COMPONENTS CONTAINING NO CRITICAL POINTS OF
A POLYNOMIAL EXCEPT FOR ITS ZEROS
\end{center}

\textit{Let $P$ be a complex polynomial of degree $n$ and let $E$
be a connected component of the set $\{z:|P(z)|\leq 1$ containing
no critical points of P different from its zeros. We prove the
inequality $|(z-a)P'(z)/P(z)|\leq n$ for all $z\in E
\backslash\{a\}$, where a is the zero of the polynomial $P$ lying
in $E$. Equality is attained for $P(z) = cz^n$ and any $z, c\neq
0$. Bibliography: 4 titles.}

\begin{center}
{\bf\Large Introduction}
\end{center}

Let $R$ be a rational function of degree $n$ represented as
$R=P/Q$, where $P$ and $Q$ are polynomials of degrees $n$ and not
exceeding $n$, respectively, which have no zeros in common. Assume
that $R(0)=R'(0)\neq 0$. Sheil-Small [1, 10.3.2] posed a question
on finding a neighborhood of the point $w=0$ in which one could
distinguish a one-valued branch, $f(w)$, of the inverse function
$z=R^{-1}(w)$, f(0)=0, satisfying the inequality
$$
{\rm Re}\,\frac{wf'(w)}{f(w)}\geq \frac1n
$$
The inequality obtained would be of interest in connection with
the shape of the level curves $|R(z)|=const$ in the domain where
the function $R$ is univalent. In the present note, a closely
related problem for polynomials $P(Q\equiv1)$ is considered. More
precisely, the following result is proved.

\textbf{Theorem.} \textit{Let $P$ be a polynomial of degree not
exceeding $n$ and let $E$ be a connected component of the
lemniscate $|P(z)|\leq 1$ containing no critical points of the
polynomial $P$ different from its zeros. Then, for any point $z\in
E\backslash\{a\}$},
$$
\left|\frac{(z-a)P'(z)}{P(z)}\right|\leq n,\eqno(1)
$$
\textit{where $a$ is the zero of the polynomial $P$ belonging to
yhe component $E$. Equality in (1) is attained for any point $z$
in the case where $P(z)=cz^n, c\neq 0$.}

If the component $E$ contains no critical points, then from (1) it
follows that for the corresponding branch $f$ of the function
inverse to the polynomial $P$, the inequality
$$
\left|\frac{wf'(w)}{f(w)-a}\right|\geq\frac1n.
$$
holds in the disk $|\omega|<1$. This result is weaker than the
Sheil-Small statement. However, contrary to [1], critical points
in $E$ are allowed. The result obtained has the following
geometric interpretation. Assume, for simplicity, that $a=0$ and
let log ($\cdot$) denote the one-valued branch of the logarithm
mapping the plane slit along the real positive semiaxis onto a
strip of width 2$\pi$. For any level curve
$c(\tau)(|P(z)|=\tau<1)$, the "curve"  $\gamma(\tau)=\log c(\tau)$
connects the opposite sides of the strip mentioned;
consequently,its length is no less than $2\pi$. On the other hand,
the image $\gamma(\tau)$ under the map $\log P(\exp(\cdot))$
covers a vertical interval of length $2\pi$ no more than $n$
times. Therefore, on the curve $\gamma(\tau)$, there is a point
$\zeta$ at which the distortion coefficient satisfies the
condition.
$$
\left|[\log P(\exp \zeta)]'\right|\leq n.
$$
This means that inequality (1) holds at a certain point
$z=\exp\zeta$ of the level curve $c(\tau)$. The theorem of the
present paper claims that this inequality holds at any point of
the curve $c(\tau)$. The assumption that the set $E$ contains no
critical points different from the polynomial zero is essential.
For instance, for the polynomial $P(z)=z^3/2-3z^2/4$, the
lemniscate $|P(z)|\leq 1$ contains both critical points $z=0$ and
$z=1$, whence it is connected. The point $z=2$ belongs to this
lemniscate, but
$$
\frac{2P'(2)}{P(2)}=6>3.
$$
\smallskip

\textbf{Corollary.} \textit{If, under the assumption of the
theorem, the inequality $P(z)>0$ holds at point $z\in E$, then the
polar derivative with respect to the point $a$, $D_aP$, satisfies
the bound}
$$
{\rm Re}\,D_aP(z)={\rm Re}\,[nP(z)-(z-a)P'(z)]\geq 0,\qquad {\rm
Im}.
$$
\smallskip

The theorem will be proved in Sec. 2. Ideologically, it comes back
to the proof of Hayman's conjecture on coverings of vertical under
a conformal mapping of the disk [2].

\begin{center}
{\bf\Large \S1. Auxiliary constructions and assertions}
\end{center}

Let $P$ be a polynomial of degree $n$ and let $E$ be a connected
component of the lemniscate $|P(z)|\leq1$ that contains no
critical points of the polynomial $P$ other than its zeros (i.e.,
no points $\zeta$ such that $P'(\zeta)=0$ and $P(\zeta)\neq0$).
Let $a$ be the zero of $P$ lying in $E$ and let $z_0$ be a point
of the component $E$ such that $P(z_0)>0$.

By $\mathscr R$ denote the Riemann surface of the function
$\mathscr P^{-1}$ inverse to the polynomial $P$. In what follws,
we consider the function $\mathscr P^{-1}$ as a one-valued
function given on the surfase $\mathscr R$. Let $\mathscr
P:\overline{\mathbb{C}}_z\to\mathscr R$ be the function inverse to
$\mathscr P^{-1}$ in this sense. The projection of a point
$W\in\mathscr R$ is defined as the point $P(\mathscr
P^{-1}(W))\in\overline{\mathbb{C}}_w$. Assume that the ray
$\{w:{\rm Im}\,w=0,\,0<{\rm Re}\,w<\infty\}$ contains no critical
values of the polynomial $P$ (i.e., no points $P(\zeta)$ such that
$P'(\zeta)=0$ for a certain $\zeta$). By $\mathscr L$ denote the
ray on the surface $\mathscr R$ or, more exactly, the Jordan curve
univalently lying over the above ray of the sphere
$\overline{\mathbb{C}}_w$ and connecting the points $\mathscr
P(a)$ and $\mathscr P(\infty)$. Let
$T=\{t_k\}_{k=0}^m,\;0=t_0<P(z_0)=t_1<\ldots<t_{m-1}<t_m=\infty$,
be a partition of the interval $0\leq  t\leq\infty$ containing all
those values of $t$ in $1<t<\infty$ at which the circle
$\gamma(t):=\{w:|w|=t\}$ contains at least one critical value
$P(\zeta)$ with $\zeta\in E$. Finally, by $\mathscr C(t)$ denote
the closed Jordan curve on $\mathscr R$ intersecting the ray
$\mathscr L$ and lying over the circle $\gamma(t)$ whose
orientation corresponds to the positive orientstion on the
projection $\gamma(t)$, $0<t<\infty,\;t\not\in T;\;c(t)$ is the
image of the curve $\mathscr C(t)$ under the mapping $\mathscr
P^{-1}$.
\smallskip

{\bf Lemma 1.} {\it The argument increment}
$$
\triangle_{ c(t)}\arg P(z)
$$
\textit{is a nondecreasing function of $t$ on the set $\{t:
0<t<\infty,\,t\not\in T\}$.}\smallskip

\textit{Proof.} Let $0<t'<t''<\infty,\;t',t''\not\in T$. The
points $\mathscr P(a)$ and $\mathscr P(\infty)$ are located on
different sides of the curves $\mathscr C(t')$ and $\mathscr
C(t'')$ on the surfase $\mathscr R$. Therfore, the nonintersecting
Jordan curves $ c(t')$ and $ c(t'')$ separate the point $a$ from
$\infty$ on the sphere $\overline{\mathbb{C}}_z$. Consequently,
one of them lies in the interior of the other. Furthermore, as we
move along the ray $\mathscr L$ from the point $\mathscr P(a)$ to
the point $\mathscr P(\infty)$, we first meet the curve $\mathscr
C(t')$ and then the curve $\mathscr C(t'')$. This means thst the
curve $c(t')$ lies in the interior of the urve $c(t'')$.
Therefore, the number $N_{t'}$ of the zeros of the polynomial $P$
lying inside $c(t')$ does not exceed the number $N_{t''}$ of the
zeros lying inside $c(t'')$ (with account for their
multiplicities). It remains to apply the argument principle:
$$
2\pi N_t=\triangle_{ c(t)}\arg P(z),\qquad t=t',t''.
$$
The lemma is proved.
\smallskip

The points $\mathscr P(a)$ and $\mathscr P(\infty)$ lie on
different sides of the curve $\mathscr C(t)$ for any $t,
0<t<\infty, t\not\in T$. This implies that for every $k=0,\ldots,
m-1$, the doubly-connected domain
$$
\mathscr G_k=\bigcup\limits_{t_k<t<t_{k+1}}\mathscr C(t)
$$
also separates the points $\mathscr P(a)$ and $\mathscr
P(\infty)$. At the same time, the curve $\mathscr P(H),\;
H=\{z:z_0+(a-z_0)\tau,\,1\leq\tau\leq\infty\}$, connects these
points. Therefore, for every $k=0,\ldots, m-1$, there is at least
one Jordane arc $\mathscr H_k$, on the curve $\mathscr P(H)$ that
lies in the domain $\mathscr G_k$ and connects its boundary
components. Thus, in the above notation, the following assertion
holds.
\smallskip

{\bf Lemma 2.} {\it For any $k=0,\ldots,m-1$ the domain $\mathscr
G_k\setminus \mathscr H_k$  is simply connected.}
\smallskip

Below, we will need the notio of condenser capacity (e.g., see
[3]). For sufficiently small positive $r$ and $\rho$ on the sphere
$\overline{\mathbb{C}}_z$, consider the condensers
$$
C(r)=(H,\{z:|z-z_0|\leq r\})
$$\hskip1cm
and
$$
C(r,\rho)=(H\cup\{z:|z-a|\leq\rho\}\cup\{z:|z|\geq
1/\rho\}\cup\bigcup\limits_{P'(\zeta)=0}\{z:|z-\zeta|\leq\rho\},
\{z:|z-1|\leq r\}).
$$
\smallskip

{\bf Lemma 3.} {\it For a fixed $r,\;0<r<|a-z_0|,$ the condensers
capacities satialy the relation}
$$
\lim\limits_{\rho\to 0}{\rm cap}\,C(r,\rho)={\rm cap}\,C(r).
$$
\smallskip

{\it Proof.} We make use of the continuity of the capacity and of
the fact that the latter is invariant under addition of a finite
number of points to points to the condenser's plates:
$$
\lim\limits_{\rho\to 0}{\rm cap}\,C(r,\rho)={\rm
cap}\,\left(H\cup\bigcup\limits_{P'(\zeta)=0}
\{\zeta\},\{z:|z-z_0|\leq r\}\right)={\rm cap}\,C(r)
$$
(see Propositions 1.4 and 1.6 in [3]). This proves the lemma.
\smallskip

Below, we introduce new notation and give some comments.

$\zeta=f_k(W)$ is the one-valued branch of the function
$\zeta=\log (W/ P(z_0))$ that maps the domain $\mathscr
G_k\setminus \mathscr H_k$ conformally and univalently into the
"strip"  $\Pi_k:=\{\zeta:\xi_{k}<{\rm
Re}\,\zeta<\xi_{k+1}\},\;k=0,\ldots,m-1$. Here, $\xi_k=\log
(t_k/P(z_0)),\;k=0,1,\ldots,m$. The choice of such a branch is
feasible in view of Lemma 2. For $k=1$ and $k=m$, $\Pi_k$ is a
half-plane.

$u(z)$ is the potential function of the condenser $C(r,\rho)$,
i.e., the resl-valued function continuous on
$\overline{\mathbb{C}}_z$, vanishing on the first plate of the
condenser $C(r,\rho)$, equal to unity on its second plate, and
harmonic in the complement to these plates:
$$
v_k(\zeta)=\left\{\begin{array}{ll} u\left(\mathscr
P^{-1}(f_k^{-1}(\zeta))\right),\qquad &\zeta\in f_k(\mathscr
G_k\setminus \mathscr H_k),\\
0,\qquad & \zeta\in \Pi_k\setminus f_k(\mathscr G_k\setminus
\mathscr H_k),\qquad k=0,\ldots,m-1\end{array}\right.
$$
On $\partial\Pi_k$ the function $v_k$ is defined by continuity.
The function obtained in this way is also denoted by $v_k$. As is
not difficult to see, the function $v_k$ satisfies the Lipschitz
condition in the strip $\overline{\Pi}_k,\;k=0,\ldots,m-1$,
whereas the function $v_j$ is equal to unity on the set
$f_j(\mathscr P(\{z:|z-z_0|\leq r\})\cap\mathscr G_j ),\;j=0,1$.

$v_k^*(\zeta)$ is the result of Steiner symmetrization of the
function $v_k(\zeta),\;\zeta\in\overline{\Pi}_k$, with respect to
the real axis (see [4]). Every function $v_k^*(\zeta)$ is a
Lipschitz function in $\overline{\Pi}_k$ and vanishes on the set
$\{\zeta\in\overline{\Pi}_k:|{\rm Im}\,\zeta|\geq \pi
n\},\;k=0,\ldots,m-1$. Lemma 1 implies the following inequalities:
$$
v_{k-1}^*(\xi_k+i\eta)\leq v_k^*(\xi_k+i\eta),\qquad
-\infty<\eta<\infty,\quad k=2,\ldots,m-1.\eqno(2)
$$

$\zeta=F(z)$ is the function that maps the unit disk $|z|<1$,
conformally and univalently, onto the strip $|{\rm Im}\,\zeta|<\pi
n$ in such a way that $F(0)=0,\;F'(0)>0$.

$\widetilde{r}$ is the upper bound for all $r$ for which the set
$F(\{z:|z|<r\})\cap \{\zeta:{\rm Re}\,(-1)^j\zeta<0\}$ belongs to
the result of Steiner symmetrization with respect to the real axis
of the set $f_j(\mathscr P(\{z:|z-z_0|\leq r\})\cap \mathscr G_j)$
for $j=0$ and $j=1$.

$v(\zeta)$ is the potential function of the condenser
$\widetilde{C}(\widetilde{r})=(\overline{\mathbb{C}}_{\zeta}\setminus
\{\zeta:|{\rm Im}\,\zeta|<\pi n\},F(\{z:|z|\leq\widetilde{r}\}))$.
It is readily seen that
$$
\begin{array}{lll}
\dfrac{\partial v}{\partial\xi}=0\quad&\mbox{on the line}\quad
{\rm
Re}\,\zeta=0,\\
\;&\;\\
 \dfrac{\partial v}{\partial\xi}\leq 0\quad&\mbox{on every line}\quad{\rm Re}\,\zeta=\xi>0.\end{array}\eqno(3)
$$
The level curves of the potential function $v$ coincide with the
level curves of the function $F$ (i.e., with the curves
$|F^{-1}(\zeta)|=\textrm{const}$).

Given a sufficiently smooth function $\lambda$ on an open set
$\Omega\subset\mathbb{C}$, denote
$$
I(\lambda,\Omega)=\int\limits_{\Omega}|\nabla\lambda|^2d\sigma.
$$
\smallskip

{\bf Lemma 4.} {\it The following inequality holds:}
$$
\sum\limits_{k=0}^{m-1} I(v_k^*,\Pi_k)\geq I(v,\mathbb{C}).
$$
\smallskip

{\it Proof.} Set $G_k=\{\zeta\in \Pi_k:|{\rm Im}\,\zeta|<\pi
n\},\;k=0,1,\ldots,m-1$, and $l_k=\{\zeta:{\rm
Re}\,\zeta=\xi_k,\;|{\rm Im}\,\zeta|<\pi n\},\;k=2,\ldots,m-1$.
For every $k,\;0\leq k\leq m-1$, we have
$$
I(v_k^*,\Pi_k)=I(v_k^*,G_k)=I(v_k^*-v+v,G_k)=I(v_k^*-v,G_k)+I(v,G_k)+
$$
$$
+2\iint\limits_{G_k}\left[\frac{\partial(v_k^*-v)}{\partial\xi}\frac{\partial
v}{\partial\xi}+\frac{\partial(v_k^*-v)}{\partial\eta}\frac{\partial
v}{\partial\eta}\right]d\xi d\eta\geq
I(v,G_k)-2\int\limits_{\partial G_k}(v_k^*-v)\frac{\partial
v}{\partial n}ds,
$$
where $\partial/\partial n$ means differentiation along the inward
normal to the boundary of the domain $G_k$ (angle points are
excluded). With account for relations (2) and (3), we derive
$$
\sum\limits_{k=0}^{m-1}I(v_k^*,\Pi_k)\geq\sum\limits_{k=0}^{m-1}
I(v,G_k)-2\sum\limits_{k=1}^{m-1}\int\limits_{\partial
G_k}(v_k^*-v)\frac{\partial v}{\partial n}ds=
$$
$$
=\sum\limits_{k=0}^{m-1}
I(v,G_k)-2\sum\limits_{k=2}^{m-1}\int\limits_{l_k}\left[(v_{k-1}^*-v)\left(-\frac{\partial
v}{\partial \xi}\right)+(v_k^*-v)\frac{\partial
v}{\partial\xi}\right] ds=
$$
$$
=\sum\limits_{k=0}^{m-1}
I(v,G_k)+2\sum\limits_{k=2}^{m-1}\int\limits_{l_k}(v_{k-1}^*-v_k^*)\frac{\partial
v}{\partial\xi}ds\geq I(v,\mathbb{C})
$$
This completes the proof.

\begin{center}
{\bf\Large \S2. Proof of the theorem}
\end{center}

It is sufficient to prove inequality (1) at an arbitrary point
$z_0\in E\setminus\{a\}$ such that $P(z_0)>0$ under the assumption
that the ray $\{w:{\rm Im}\,w=0,\,0<{\rm Re}\,w<\infty\}$ contains
no critical points of the polynomial $P$. We use the notation
introduced in Sec. 1. The chain of relations below stems from the
conformal invariance of the Dirichlet integral, from the P\'{o}lya
and Szeg\H{o} theorem on function symmetrization (see [4]), and
from Lemma 4:
$$
{\rm cap}\,C(r,\rho)=I(u,\mathbb{C})\geq\sum\limits_{k=0}^{m-1}
I(v_k,\Pi_k)\geq\sum\limits_{k=0}^{m-1} I(v_k^*,\Pi_k)\geq
$$
$$
\geq I(v,\mathbb{C})={\rm cap}\,\widetilde{C}(\widetilde{r}).
$$
In view of Lemma 3, we ultimately obtain
$$
{\rm cap}\,C(r)\geq{\rm
cap}\,\widetilde{C}(\widetilde{r}).\eqno(4)
$$
In order to compute the asymptotics of the condenser capacity as
$r\to 0$, we use known formulas (e.g., see [3, (1.6) and (1.8)]),
in which $r(B,a)$ stands for the inner radius of the domain $B$
with respect to a point $a\in B$. As a result, we obtain
$${\rm cap}\,C(r)=-\frac{2\pi}{\log r}-\frac{1}{2\pi}(\log
r(\mathbb{C}_z\setminus H,z_0))\left(\frac{2\pi}{\log
r}\right)^2+o\left(\left(\frac{1}{\log r}\right)^2\right)=
$$
$$
=-\frac{2\pi}{\log r}-2\pi(\log[4|a-z_0|])\left(\frac{1}{\log
r}\right)^2+o\left(\left(\frac{1}{\log r}\right)^2\right),\qquad
r\to 0.
$$
Further, the second plate of the condenser
$\widetilde{C}(\widetilde{r})$ is an "almost disk" of radius
$(r|P'(z_0)|/P(z_0))(1+o(1))$ as $r\to 0$. Consequently,
$$
{\rm
cap}\,\widetilde{C}(\widetilde{r})=-\frac{2\pi}{\log(r|P'(z_0)|/P(z_0))}-
$$
$$
-2\pi(\log r(\{\zeta:|{\rm Im}\,\zeta|<\pi
n\},0))\left(\frac{1}{\log (r|P'(z_0)|/P(z_0))}\right)^2+
$$
$$
+o\left(\left(\frac{1}{\log r}\right)^2\right)=-\frac{2\pi}{\log
(r|P'(z_0)|/P(z_0))}-
$$
$$
-2\pi(\log(4n))\left(\frac{1}{\log(r|P'(z_0)|/P(z_0))}\right)^2+o\left(\left(\frac{1}{\log
r}\right)^2\right)=
$$
$$
=-\frac{2\pi}{\log r}\left(1-\frac{\log|P'(z_0)/P(z_0)|}{\log
r}+o\left(\frac{1}{\log r}\right)\right)-
$$
$$
-2\pi(\log(4n))\left(\frac{1}{\log
r}\right)^2+o\left(\left(\frac{1}{\log
r}\right)^2\right)=-\frac{2\pi}{\log r}-
$$
$$
-2\pi(\log|4nP(z_0)/P'(z_0)|)\left(\frac{1}{\log
r}\right)^2+o\left(\left(\frac{1}{\log r}\right)^2\right),\qquad
r\to 0.
$$
Substituting the asymptotics obtained into inequality (4), we
arrive at the inequality
$$
|a-z_0|\leq|nP(z_0)/P'(z_0)|.
$$
Under the assumptions considered, the latter relation coincides
with (1) $(z=z_0)$. The case of equality is verified
straightforwardly.

The theorem is proved.

\begin{center}
{\bf\Large Литература}
\end{center}

\begin{enumerate}
\item T.Sheil-Small, {\it Complex polynomials}, Cambridge Univ.
Press., Cambridge (2002).
\item V.N. Dubinin, "Coverings of vertical segments under a conformal mapping", \textit{Mat. Zametki} {\bf 28}, No. 1, 25-32 (1980).
\item V.N. Dubinin, "Symmetrization in the geometric theory of function of a complex variable", \textit{Usp. Mat. Nauk}, {\bf 49}, No. 1,
3-76 (1994).
\item W.K.Hayman, {\it Multivalent functions}, Cambridge Univ.
Press., Cambridge (1994).
\end{enumerate}

\end{document}